\font\Bbb=msbm10 
\def\matC{\hbox{\Bbb C}}

\font\Gros=cmbx10 scaled\magstep1

\magnification \magstep1 

\null  
\bigskip
\bigskip
\centerline{\Gros Sur la caract\'erisation des } 
\smallskip
\centerline{\Gros r\'etractes holomorphes \`a l'aide} 
\smallskip

\centerline{\Gros de la m\'etrique infinit\'esimale de Kobayashi}  

\medskip
\centerline{\Gros Jean-Pierre Vigu\'e} 
\bigskip

{\bf1. Introduction} 
\medskip

Dans un article r\'ecent avec M. Abate [2], nous avons montr\'e le
 th\'eor\`eme suivant. 
\medskip

{\bf Th\'eor\`eme 1.1.}{\it- Soient \/}$ (E_{1},\|.\|_{1}) $ {\it et \/}$ (E_{2},\|.\|_{2})
 $ {\it deux espaces de Banach complexes et soient \/}$ B_{1} $ {\it
 et \/}$ B_{2} $ {\it leurs boules-unit\'es ouvertes. 
Soit \/}$ f:B_{1}\longrightarrow B_{2} $ {\it une application holomorphe telle
 que \/}$ f(0)=0 $ {\it et que \/}$ f' (0) $ {\it soit une isom\'etrie
 (c'est-\`a-dire que, pour tout \/}$ x\in E_{1}${\it, \/}$\|f'(0).x\|_{2}=\|x\|_{1}${\it). Alors,
 les assertions suivantes sont \'equivalentes :\/} 

(i) {\it il existe une d\'ecomposition directe \/}$ E_{2}=f'(0)(E_{1})\oplus
 F $ {\it telle que la projection associ\'ee \/}$\pi:E_{2}\longrightarrow
 f'(0)(E_{1}) $ {\it soit de norme \/}$ 1${\it,\/} 

(ii) $ f(B_{1}) $ {\it est une sous\/}$-${\it vari\'et\'e directe
 ferm\'ee de \/}$ B_{2}${\it, \/}$ f $ {\it est un biholomorphisme de
 \/}$ B_{1} $ {\it sur \/}$ f(B_{1}) $ {\it et il existe une r\'etraction
 holomorphe de \/}$ B_{2} $ {\it sur \/}$ f(B_{1})${\it.\/} 
\medskip

Signalons aussi le r\'esultat suivant, d\^u \`a M. Abate [1] : Soit
 $ D $ un domaine born\'e taut de $ {\matC}^{n} $ et 
soit $ f:D\longrightarrow D $ une application holomorphe. Consid\'erons la
 suite des it\'er\'ees $ f^{n} $ de $ f$. Alors, si la 
suite $f^n$ n'est pas compactement divergente, il existe une r\'etraction
 holomorphe unique $ \rho:D\longrightarrow D $ adh\'erente
 \`a la suite $ f^{n} $ et $ f $ est un automorphisme analytique de $ \rho(D)$.

\medskip

Ces deux r\'esultats montrent l'importance des r\'etractions holomorphes.
 Dans cet article, nous allons donner une caract\'erisation des sous-vari\'et\'es
 complexes ferm\'ees $ V $ de la boule-unit\'e ouverte de $ {\matC}^{n}
 $ pour une norme, r\'etractes holomorphes 
de $D$ en utilisant la m\'etrique infinit\'esimale de Kobayashi.
 Nous donnerons aussi un certain nombre d'exemples et d'applications.
\medskip

{\bf2. Rappels} 
\medskip
On dit qu'une application holomorphe 
$ \rho:D\longrightarrow D $ 
est une r\'etraction holomorphe si $ \rho^{2}=\rho$, ou, ce qui revient
 au m\^eme que $ \rho $ est \'egal \`a l'identit\'e sur son image $ \rho(D)$.
 On dit alors que son image $ \rho(D)$, qui est une sous-vari\'et\'e
 analytique complexe de $ D$, est r\'etracte holomorphe de $ D$.
\medskip

Pour la d\'efinition et les propri\'et\'es des m\'etriques infinit\'esimales
 de Carath\'eodory $ E_{D} $ et de Kobayashi $ F_{D} $ sur un domaine born\'e
 $ D $ de $ {\matC}^{n} $ (ou plus g\'en\'eralement, une vari\'et\'e analytique
 complexe), nous renvoyons le lecteur au livre de S. Kobayashi [6]
 (on peut aussi consulter les livres de T. Franzoni et E. Vesentini
 [4] ou de M. Jarnicki et P. Pflug [5]). 
\medskip

On sait d'autre part que, si $ f $ est une application holomorphe de
 $ D_{1} $ dans $ D_{2}$, on a, pour tout $ x\in D_{1}$, pour tout $ v $ appartenant
 \`a l'espace tangent $ T_{x}(D_{1})$, 
$$F_{D_{2}}(f(x),f'(x).v)\leq F_{D_{1}}(x,v)\hbox{\rm,}$$
 $$E_{D_{2}}(f(x),f'(x).v)\leq E_{D_{1}}(x,v)\hbox{\rm.}$$
 En particulier, si $ f $ est un isomorphisme analytique de $ D_{1} $ sur
 $ D_{2}$, alors les in\'egalit\'es pr\'ec\'e\-den\-tes sont des \'egalit\'es.
 
\medskip

Enfin, remarquons que, si $ B $ est la boule-unit\'e ouverte de $ {\matC}^{n} $ 
pour une norme $ \|.\|$, on d\'eduit du th\'eor\`eme de Hahn-Banach
 que, pour tout $ v\in{\matC}^{n}$, 
$$E_{B}(0,v)=F_{B}(0,v)=\|v\|\hbox{\rm.}$$
\medskip

 Rappelons enfin le r\'esultat de L. Lempert [7 et 8] (voir aussi
 le livre de S. Kobayashi [6]) : soit $ D $ un domaine born\'e convexe
 de $ {\matC}^{n}$. Sur $ D$, les m\'etriques infinit\'esimales de Carath\'eodory
 $ E_{D} $ et de Kobayashi $ F_{D} $ co\"{\i}ncident. 
\medskip

{\bf3. Les r\'etractes holomorphes de dimension 1} 
\medskip

Nous allons commencer par donner les r\'esultats en dimension 1. Comme
 cela appara\^{\i}tra clairement dans la suite (voir aussi L. Lempert
 [8]), c'est un cas fondamentalement diff\'erent de celui de la dimension
 sup\'erieure et, du moins, dans le cas o\`u $ D $ est un domaine born\'e convexe,
 il est assez bien compris. Cependant, il me semble int\'eressant de reprendre les
 r\'esultats obtenus en dimension 1 [9]. 
\medskip
Rappelons qu'une g\'eod\'esique complexe (au sens de E. Vesentini)
 d'un domaine born\'e $D$ de $ {\matC}^{n} $ est une application holomorphe
 $ \varphi:\Delta\longrightarrow D $ qui est une isom\'etrie
 pour les m\'etriques infinit\'esimales de Carath\'eodory en un point de 
$ \Delta $ (ou en tout point de $ \Delta $ ; les deux 
propri\'et\'es sont \'equivalentes).
\medskip

Consid\'erons un domaine born\'e convexe $ D $ de $ {\matC}^{n}$. Soit
 $ V $ une sous-vari\'et\'e complexe de $ D$. Si $ V $ est r\'etracte holomorphe
 de $ D$, alors il est facile de voir que $ V $ est simplement connexe.
 Si on suppose de plus que $ V $ est de dimension 1, ceci entra\^{\i}ne
 que $ V $ est analytiquement isomorphe au disque-unit\'e $ \Delta $ et
 que $ V $ est l'image d'une g\'eod\'esique complexe 
$ \varphi:\Delta\longrightarrow D $ au sens de E. Vesentini. R\'eciproquement,
 si $ \varphi:\Delta\longrightarrow D $ est une g\'eod\'esique
 complexe, son image $ \varphi(\Delta) $ est isomorphe au disque-unit\'e
 $ \Delta $ et est r\'etracte holomorphe de $ D$. Dans ce cas-l\`a, on
 peut montrer que $ D $ admet beaucoup de g\'eod\'esiques complexes et
 donc de r\'etractions holomorphes, alors que, comme l'a d\'ej\`a remarqu\'e
 L. Lempert [8], en dimensions strictement sup\'erieures \`a 1, les
 r\'etractions holomorphes n'ont pas de raison d'exister. 
\medskip

Si on veut essayer de caract\'eriser les r\'etractions holomophes
 de dimension 1 dans des domaines born\'es non n\'ecessairement convexes,
 il peut \^etre utile de d\'efinir d'autres m\'etriques invariantes
 [9], mais, pour l'instant, les r\'esultats obtenus sont techniques
 et compliqu\'es. 
\medskip

{\bf4. Le r\'esultat principal et sa d\'emonstration} 
\medskip

Le r\'esultat principal que nous allons d\'emontrer dans cet article
 est le suivant. 
\medskip

{\bf Th\'eor\`eme 4.1.}{\it- Soit \/}$ B $ {\it la boule-unit\'e ouverte
 de \/}${\matC}^{n}${\it, muni d'une norme \/}$\|.\|${\it, et soit
 \/}$ V $ {\it une sous-vari\'et\'e complexe ferm\'ee de \/}$ B $ {\it
 telle que \/}$ 0\in V${\it. Pour que \/}$ V $ {\it soit r\'etracte holomorphe
 de \/}$ V${\it, il faut et il suffit que les deux conditions suivantes
 soient satisfaites :\/}

(i) $ F_{B}(0,v)=F_{V}(0,v), \forall v\in T_{0}(V),$ 

(ii) {\it il existe un projecteur \/}$ p $ {\it de norme \/}$ 1 $ {\it
 de \/}${\matC}^{n} $ {\it sur \/}$ T_{0}(V) $ {\it(\/}${\matC}^{n} $ {\it
 et \/}$ T_{0}(V)$\'e{\it tant munis de la norme \/}$\|.\|${\it).\/}
 \medskip

Avant de d\'emontrer ce r\'esultat, remarquons que, d'apr\`es le r\'esultat
 que nous avons rappel\'e, on sait que $ F_{B}(0,.)=\|.\|$. Par suite
 (i) signifie que $ F_{V}(0,.) $ est \'egal \`a la norme $ \|.\| $ restreinte
 \`a $ T_{0}(V)$. 
\medskip

{\it D\'emonstration\/}. Supposons d'abord qu'il existe une r\'etraction
 holomorphe $ \rho: B\longrightarrow V$. En consid\'erant
 l'injection $ i:V\longrightarrow D$, on obtient : 
$$F_{V}(0,(\rho{\circ}i)'(0).v)\leq F_{B}(0,i'(0).v)\leq F_{V}(0,v)\hbox{\rm.}$$
 Mais $ (\rho{\circ}i)' (0)={\rm id}$. Par suite, pour tout $ v \in T_{0}(V)$,
 on a : $ F_{B}(0,v)=F_{V}(0,v)$, ce qui montre (i). 
\medskip

D'apr\`es les propri\'et\'es de la m\'etrique infinit\'esimale de
 Kobayashi, il est clair que $ \rho' (0) $ est un projecteur de norme
 $ 1 $ de $ {\matC}^{n} $ sur $ T_{0}(V)$, quand on munit $ {\matC}^{n}
 $ de la m\'etrique infinit\'esimale de Kobayashi $ F_{B}(0,.) $ et $ T_{0}(V)
 $ de la m\'etrique infinit\'esimale de Kobayashi $ F_{{\rm V}}(0,.)$.
 Mais, nous venons de montrer que $ F_{B}(0,.) $ est \'egal \`a la norme
 $ \|.\| $ et que $ F_{{\rm V}}(0,.) $ est la restriction \`a $ T_{0}(V)
 $ de la norme $ \|.\|$. Ceci suffit \`a d\'emontrer que $ \rho' (0) $ est
 un projecteur de norme $ 1 $ pour la norme $ \|.\|$. 
\medskip

Montrons maintenant la r\'eciproque. Remarquons d'abord que le projecteur
 $ p $ de norme $ 1 $ dont nous avons suppos\'e l'existence envoie $ B $ dans
 $ B$. Aussi, $ p|_{B}:B\longrightarrow B $ est une r\'etraction
 holomorphe de $ B $ sur $ B\cap p({\matC}^{n})=B\cap T_{0}(V)$. Consid\'erons
 maintenant $ p|_{V}:V\longrightarrow B\cap T_{0}(V)$.
 Comme $ B\cap T_{0}(V) $ est la boule-unit\'e ouverte de $ T_{0}(V)$,
 ceci entra\^{\i}ne que $ F_{B\cap T_{0}(V)}(0,.)=\|.\|$. On en d\'eduit
 que, pour tout $ v\in T_{0}(V)$, $ F_{V}(0,v)=F_{B\cap T_{0}(V)}(0,v)$.
\medskip

Aussi, $ p|_{V}:V\longrightarrow B\cap T_{0}(V) $ est une
 isom\'etrie pour la m\'etrique infinit\'esimale de Kobayashi \`a l'origine,
 ce qui signifie que$, $ pour tout $ v\in T_{0}(V)$, 
$$F_{B\cap T_{0}(V)}(0,v)=F_{V}(0,v)=F_{V}(0,p(v)).$$
 D'autre part, on sait que $ V $ qui est une sous-vari\'et\'e complexe
 ferm\'ee de $ B $ est compl\`ete pour la distance de Kobayashi et est
 donc taut. On peut donc appliquer le th\'eor\`eme de L. Belkhchicha
 [3]. Ainsi, $ p|_{V} $ est un isomorphisme analytique de $ V $ sur $ B\cap
 T_{0}(V)$. Alors l'application $ \rho=(p|_{V})^{-1}{\circ}p $ est bien
 une r\'etraction holomorphe de $ B $ sur $ V$, et le th\'eor\`eme est
 d\'emontr\'e. 
\medskip
Ce r\'esultat peut bien s\^ur s'appliquer quand $ V $ est de dimension
 1. Il peut \^etre int\'eressant de le traduire dans ce cas. On d\'eduit
 facilement du th\'eor\`eme de Hahn-Banach que, si $ F $ est un sous-espace
 vectoriel de dimension 1 de $ {\matC}^{n}$, il existe une projection
 $ p:{\matC}^{n}\longrightarrow F $ de norme 1. La condition
 (i) du th\'eor\`eme 4.1 (portant sur l'\'egalit\'e des m\'etriques
 infinit\'esimales de Kobayashi) suffit donc \`a entra\^{\i}ner que $ V $ est
 r\'etracte holomorphe de $ B $ (et aussi que $ V $ est l'image d'une g\'eod\'esique
 complexe de $ B$).
\medskip

{\bf5. Applications} 
\medskip

A titre d'applications, nous allons traiter quelques exemples. 
\medskip

Exemple 5.1.- Soit $ B $ la boule-unit\'e ouverte de $ {\matC}^{2}$,
 muni de la norme hermitienne et soit 
$$V=B\cap\{(x,y)\in{\matC}^{2} | y=x^{2}\}\hbox{\rm.}$$
 Alors $ V $ n'est pas r\'etracte holomorphe de $ B$. 

En effet, $ V $ est l'ensemble des $ (x,y)\in{\matC}^{2}$ tels que $ y=x^{2}
 $ et que $ |x|^{2}+|x|^{4}<1 $ ou, si l'on pr\'ef\`ere que $ |x|$ soit
 strictement inf\'erieur \`a une constante $ \rho<1$. Il est alors facile
 de v\'erifier que, pour tout $ v\in T_{0}(V), (v\neq0)$, $ F_{V}(0,v)\neq
 F_{B}(0,v)$. Ainsi, $ V $ n'est pas r\'etracte holomorphe de $ B$.  
\medskip

Nous allons voir maintenant comment le le th\'eor\`eme 4.1 permet
 de ramener les questions \`a des probl\`emes pour des applications
 lin\'eaires, ce qui est plus facile. 
\medskip

Exemple 5.2.- Soit $ B $ le polydisque 
ouvert $ \Delta^{3}\subset{\matC}^{3}$ et consid\'erons l'application 
lin\'eaire $ \varphi:(x,y)\mapsto (x,y,x+y)$. 
Il est clair que 
$ V=\varphi({\matC}^{2})\cap\Delta^{3} $ est une sous-vari\'et\'e complexe ferm\'ee
 de $ B$. L'ensemble $ (1,0,1)+(\{0\}\times\Delta\times\{0\}) $ est contenu
 dans la fronti\`ere de $ B$. S'il existait une projection $ p $ de norme
 1 de $ {\matC}^{3} $ sur $ \varphi({\matC}^{2})$, cette projection
 devrait s'annuler sur $ \{0\}\times{\matC}\times\{0\}$. En consid\'erant
 le point $ (0,1,1)$, on trouve un autre sous-espace de $ {\matC}^{3}$.
 Ainsi, une telle projection ne peut pas exister ! 
\medskip

En appliquant le th\'eor\`eme 4.1, on en d\'eduit que la sous-vari\'et\'e
 $ V $ n'est pas r\'etracte holomorphe de $ B$. 
\medskip

{\bf6. Questions} 
\medskip

Les r\'esultats pr\'ec\'edents am\`enent \`a se poser la question
 suivante : soit $ D $ un domaine born\'e convexe de $ {\matC}^{n}$. Soit
 $ V $ une sous-vari\'et\'e complexe ferm\'ee de $ D $ et soit $ a\in V. $ 
Il est facile de voir que, s'il existe une r\'etraction holomorphe
 $ \rho:D\longrightarrow V$, alors on a : 

(i) $ F_{V}(a,.)=F_{D}(a,.)|_{T_{a}(V)}$, 

(ii) il existe une projection $p$ de norme 1, lorsque $ {\matC}^{n} $ 
et $ T_{a}(V) $ sont munis des m\'etriques infinit\'esimales de Kobayashi
 $ F_{D}(a,.) $ et $ F_{V}(a,.)$. 
\medskip

On peut se demander si la r\'eciproque est exacte, ce qui g\'en\'eraliserait
 le th\'eor\`eme 4.1. Pour l'instant, ce r\'esultat semble difficile
 \`a aborder. il faudrait en particulier g\'en\'eraliser le th\'eor\`eme
 de L. Belkhchicha [3], ce qui est sans doute assez d\'elicat. 
\medskip

\centerline{\bf Bibliographie} 
\smallskip

1. M. Abate. Iteration theory, compactly divergent sequences and
 commuting holomorphic maps. Ann. Scuola Norm. Sup. Pisa Cl. Sci.
 (4), {\bf18 }(1991), 167--191. 
\smallskip

2. M. Abate and J.-P. Vigu\'e. Isometries for the Carath\'eodory metric.
 Proc. AMS, \`a para\^{\i}tre. 
\smallskip

3. L. Belkhchicha. Caract\'erisation des isomorphismes analytiques
 sur la boule-unit\'e de $ {\matC}^{n} $ pour une norme. Math. Z., 
{\bf215 }(1994), 129--141. 
\smallskip

4. T. Franzoni and E. Vesentini. Holomorphic maps and invariant distances.
 Notas de Matematica [Mathematical Notes], {\bf69}. North-Holland Publishing
 Co, Amsterdam, 1980. 
\smallskip

5. M. Jarnicki and P. Pflug. Invariant distances and metrics in complex
 analysis. de Gruyter Expositions in Mathematics, {\bf9}, Walter de
 Gruyter Co, Berlin, 1993. 
\smallskip

6. S. Kobayashi. Hyperbolic complex spaces, Grundlehren 
der Mathematischen Wissenschaften [Fundamental Principles 
of Mathematical Sciences], 318, Springer-Verlag, Berlin, 1998.
\smallskip

7. L. Lempert. La m\'etrique de Kobayashi et la 
repr\'esentation des domaines sur la boule.
 Bull. Soc. Math. France {\bf109} (1981), 427--474. 
\smallskip

8. L. Lempert. Holomorphic retracts and intrinsic metrics in convex domains. 
Anal. Math. {\bf 8} (1982), 257--261. 
\smallskip

9. J.-P. Vigu\'e. G\'eod\'esiques complexes et r\'etractes holomorphes de dimension
 1. Ann. Mat. Pura Appl. (4), {\bf176 }(1999), 95--112. 
\medskip

Jean-Pierre Vigu\'e 

UMR CNRS 6086 

Universit\'e de Poitiers 

Math\'ematiques 

SP2MI, BP 30179 

86962 FUTUROSCOPE  

e-mail : vigue@math.univ-poitiers.fr

ou jp.vigue@orange.fr

 \bye